\documentclass[12pt]{amsart}
\usepackage{amssymb}
\usepackage{color}
\usepackage{amsmath,epic,curves,amscd}
\usepackage[english]{babel}
\usepackage{graphicx}
\usepackage{comment}
\usepackage{appendix}
\usepackage{mathdots}
\usepackage[all]{xy}
\usepackage{mathtools}
\usepackage{lscape}
\usepackage{stmaryrd}
\usepackage{mathabx}
\newtheorem{claim}{}[section]
\newtheorem{theorem}[claim]{Theorem}

\newtheorem{conjecture}[claim]{Conjecture}
\newtheorem{lemma}[claim]{Lemma}
\newtheorem{proposition}[claim]{Proposition}
\newtheorem{corollary}[claim]{Corollary}
\newtheorem{example}[claim]{Example}
\theoremstyle{remark}

\renewenvironment{proof}{\noindent{\it Proof. \hskip0pt}}
                      {$\square$\par\medskip}
\allowdisplaybreaks \textwidth 15.5 true cm \textheight 23.9 true cm
\usepackage{color}
\newcommand{\markg}[1]{{\color{black}#1}}

\begin{document}
\baselineskip 6.0 truemm
\parindent 1.5 true pc

\newcommand\lan{\langle}
\newcommand\ran{\rangle}
\newcommand\tr{\operatorname{Tr}}
\newcommand\ot{\otimes}
\newcommand\ttt{{\text{\rm t}}}
\newcommand\rank{\ {\text{\rm rank of}}\ }
\newcommand\choi{{\rm C}}
\newcommand\dual{\star}
\newcommand\flip{\star}
\newcommand\cp{{\mathbb C}{\mathbb P}}
\newcommand\ccp{{\mathbb C}{\mathbb C}{\mathbb P}}
\newcommand\pos{{\mathcal P}}
\newcommand\tcone{T}
\newcommand\mcone{K}
\newcommand\superpos{{\mathbb S\mathbb P}}
\newcommand\blockpos{{\mathcal B\mathcal P}}
\newcommand\jc{{\text{\rm JC}}}
\newcommand\dec{{\mathbb D}{\mathbb E}{\mathbb C}}
\newcommand\ppt{{\mathcal P}{\mathcal P}{\mathcal T}}
\newcommand\pptmap{{\mathbb P}{\mathbb P}{\mathbb T}}
\newcommand\xxxx{\bigskip\par ================================}
\newcommand\join{\vee}
\newcommand\meet{\wedge}
\newcommand\ad{\operatorname{Ad}}
\newcommand\ldual{\varolessthan}
\newcommand\rdual{\varogreaterthan}
\newcommand{\slmp}{{\mathcal M}^{\text{\rm L}}}
\newcommand{\srmp}{{\mathcal M}^{\text{\rm R}}}
\newcommand{\smp}{{\mathcal M}}

\title{Convex cones in mapping spaces between matrix algebras}

\author{Mark Girard, Seung-Hyeok Kye and Erling St\o rmer}
\address{Institute for Quantum Computing and School of Computer Science, University of Waterloo, Waterloo, ON N2L 3G1, Canada}
\email{mark.girard at uwaterloo.ca}
\address{Department of Mathematics and Institute of Mathematics, Seoul National University, Seoul 151-742, Korea}
\email{kye at snu.ac.kr}
\address{Department of Mathematics, University of Oslo, 0316 Oslo, Norway}
\email{erlings at math.uio.no} 

\thanks{SHK was partially supported by NRF-2017R1A2B4006655, Korea}

\begin{abstract}
We introduce the notion of \emph{one-sided mapping cones} of positive linear maps between matrix
algebras. These are convex cones of maps that are invariant under compositions by completely positive maps
from either the left or right side. The duals of such convex cones can be characterized
in terms of ampliation maps, which can also be used to characterize many notions from quantum information
theory---such as separability, entanglement-breaking maps, Schmidt numbers, as well as
decomposable maps and $k$-positive maps in functional analysis.
In fact, such characterizations hold if and only if the involved cone is a one-sided mapping cone.
Through this analysis, we obtain
mapping properties for compositions of cones from which we also obtain several equivalent statements
of the PPT (positive partial transpose) square conjecture.
\end{abstract}
\maketitle

\section{Introduction}

A \emph{mapping cone} is a closed convex cone of positive linear maps that is
closed under compositions by completely positive linear maps from both sides.
The notion of mapping cones was introduced by
the third author \cite{{stormer-dual}} in the 1980s to study extension problems
of positive linear maps and has been studied in the context of quantum information
theory \cite{{sko-laa},{ssz},{stormer08},{stormer08-2},{stormer08-3},{stormer_book},{stormer18}}.
Various notions from quantum information theory---such as separability, Schmidt numbers,
positive partial transposes, and entanglement-breaking maps---can be explained in terms of mapping cones.
The study of mapping cones is also  closely related to
that of operator systems~\cite{johnston_oper_sys}.

In this paper, we consider closed convex cones of positive linear maps between matrix algebras
which are closed under compositions by completely positive maps
from only one side. That is, cones $K$ that satisfy either $K\circ\cp\subset K$ or $\cp\circ K\subset K$,
respectively, where $\cp$ denotes the convex cone of all
completely positive maps and $K_1 \circ K_2$ is the set of all maps of the form
$\phi_1 \circ \phi_2$ for maps $\phi_1 \in K_1$ and $\phi_2 \in K_2$.
Mapping cones satisfy many nice properties
regarding the compositions of maps \cite{{stormer08-3},{sko-laa},{stormer_scand_2012},{stormer18}}.
We show that these properties can also be used to characterize one-sided mapping cones.

The notion of ampliation (i.e., $1\ot \phi$ or $\phi\ot 1$) of a linear map $\phi$ by the identity map on matrix algebras
plays an important role in operator algebras, as evidenced by
Stinespring's representation theorem \cite{stine} and a characterization of decomposability
due to the third author \cite{stormer82}.
It was also shown by the second author \cite{eom-kye} that $k$-positivity of linear maps can be characterized by
ampliation. Moreover, ampliation is also useful for characterizing certain kinds of linear maps in quantum information
theory---such as entanglement-breaking
maps \cite{{holevo_BE},{shor},{hsrus}}---and is crucial for distinguishing several kinds of positive (semi-definite)
matrices in tensor product spaces. For example, some criteria for separability of quantum states can be presented
in terms of ampliation of positive maps \cite{{horo-1}}.
In this paper, we provide a single framework that allows us to recover all of the above-mentioned results. The main idea is that
these kinds of characterizations hold only when the involved convex cones are one-sided mapping cones.

For this purpose, we first study the dual cones of closed convex cones of (not-necessarily positive) maps
and their compositions $K_1 \circ K_2$. We obtain various relations among dual cones and
compositions of convex cones, from which we also obtain several equivalent statements to the PPT-square conjecture.

Throughout this paper, we denote by $M_A$ the matrix algebra acting on the finite-dimensional Hilbert space
$\mathbb C^A$. We will work in the real vector spaces $H(M_A,M_B)$
consisting of all Hermitian-preserving linear maps from $M_A$
into another matrix algebra $M_B$.
Recall that a linear map $\phi:M_A\to M_B$ is called {\sl Hermitian-preserving} if
$a=a^*$ implies $\phi(a)=\phi(a)^*$---or equivalently if $\phi(a^*)=\phi(a)^*$---holds for every $a\in M_A$.

We \markg{define a bilinear pairing on the matrix algebra $M_A$ as}
\begin{equation}\label{bi-linear}
\markg{\lan a,b\ran_A=\tr (a^{\ttt} b)=\sum_{i,j}a_{ij}b_{ij}}
\end{equation}
for every $a=[a_{ij}]$ and $b=[b_{ij}]$ in $M_A$\markg{, where $a^\ttt$ denotes the transpose of $a$.
The analogous bilinear pairings on $M_B$ and $M_A\ot M_B$ will be denoted by
$\lan\ ,\ \ran_B$ and $\lan\ ,\ \ran_{AB}$, respectively. Note that this is not the Hilbert--Schmidt
inner product typically used in the literature. Using the pairing defined in~\eqref{bi-linear}
allows us to simplify many identities involving duality of cones considered in this paper.
The main results of this paper are not affected by this choice of bilinear pairing, as most
of the cones of matrices that we consider satisfy $K=K^\ttt = \overline{K}$ (e.g., the cone of
positive matrices). Moreover, this pairing is an inner product when restricted to the real
vector space of Hermitian matrices in $M_A$.}

In the next section we present relationships between the dual cones and compositions $K_1 \circ K_2$ for subsets $K_1$ and $K_2$
of linear maps that characterize one-sided mapping cones. In Section \ref{sec-mp} we revisit properties
of mapping cones in relation to one-sided mapping cones, after which we characterize the duals of
one-sided mapping cones in terms of ampliation maps in Section \ref{sec-amp}. This allows us to recover
several well-known results in quantum information theory. Finally,
in Section \ref{sec-PPT} we present several equivalent formulations of the PPT-square conjecture.

An earlier version of this paper due to second and third authors originally appeared as
arXiv:2002.09614v1. After the first author joined as a co-author, the authors changed the title, rewrote and extended
the results of the earlier version.

\section{One-sided mapping cones}

For a given linear map $\phi:M_A\to M_B$, its Choi matrix $\choi_\phi\in M_A\ot M_B$ is defined as
$$
\choi_\phi=\sum_{ij}e_{ij}\ot \phi(e_{ij})\in M_A\ot M_B,
$$
where $\{e_{ij}\}$ denote the matrix units in $M_A$. The mapping defined by $\phi\mapsto \choi_\phi$ is a linear isomorphism
between the vector space $L(M_A,M_B)$ of all linear maps from $M_A$ to $M_B$ onto the tensor product
$M_A\ot M_B$ and is usually called the Jamio\l kowski-Choi isomorphism \cite{{dePillis},{jam_72},{choi75-10}}.
It is easy to see that a map $\phi$ is Hermitian-preserving if and only if $\choi_\phi$ is a self-adjoint matrix.
Recall that $\phi$ is completely positive if and only if $\choi_\phi$ is positive (semi-definite) \cite{choi75-10}.
Note that, for all $a\in M_A$ and $b\in M_B$, we have the identity
\begin{equation}\label{dual-id}
\lan a\ot b, \choi_\phi\ran_{AB}
=\sum_{i,j}\lan b,\phi(e_{ij})\ran_B
=\lan b,\phi( a)\ran_B.
\end{equation}
We refer the reader to Chapter~4 of \cite{stormer_book} for further properties of Choi matrices.

For a given map $\phi\in L(M_A,M_B)$, \markg{we define its \emph{adjoint}} map $\phi^*:M_B\to M_A$ by the condition
\[
\lan \phi(a),b\ran_B=\lan a,\phi^*(b)\ran_A,\qquad \text{for all }a\in M_A, b\in M_B.
\]
It is easy to see that $\phi\in H(M_A,M_B)$ if and only if $\phi^*\in H(M_B,M_A)$.
For every $a\in M_A$ and $M_B$, note that
\begin{equation}\label{eq:ChoiPhiStarFlip}
\markg{\lan b\ot a,\choi_{\phi^*}\ran_{BA}
=\lan a,\phi^*(b)\ran_A=\lan b,\phi(a)\ran_B=\lan a\ot b,\choi_\phi\ran_{AB}},
\end{equation}
and thus the Choi matrix $\choi_{\phi^*}\in M_B\ot M_A$ of $\phi^*$ is the flip of \markg{$\choi_\phi$. Before proceeding,
we remark that the adjoint defined in this way is not the same as the adjoint with respect
to the Hilbert--Schmidt inner product that is commonly used in the literature with the same notation.
However, as mentioned previously, using the adjoint defined in this manner allows us to simplify many of the main results.}

We define \markg{a bilinear pairing} on $L(M_A,M_B)$ as
\begin{equation}\label{eq:maps_blinear_pair}
\lan\phi,\psi\ran=\lan\choi_\phi,\choi_\psi\ran_{AB}=\sum_{i,j} \lan \phi(e_{ij}),\psi(e_{ij})\ran_B,
\end{equation}
for every $\phi,\psi\in L(M_A,M_B)$ (see, e.g., \cite{sko-laa}). It is clear that we have the identity
\begin{equation}\label{11}
\lan\phi,\psi\ran=\lan\psi^*,\phi^*\ran.
\end{equation}
For maps $\phi\in L(M_A,M_B)$, $\psi\in L(M_B,M_C)$ and $\sigma\in L(M_A,M_C)$,
we also have that
\[
\lan\psi\circ\phi,\sigma\ran
=\sum_{i,j}\lan\psi(\phi(e_{ij})),\sigma(e_{ij})\ran_C
=\sum_{i,j}\lan \phi(e_{ij}),\psi^*(\sigma(e_{ij}))\ran_B
=\lan\phi,\psi^*\circ\sigma\ran,
\]
which further implies that $\lan\phi,\psi^*\circ\sigma\ran
=\lan \sigma^*\circ\psi,\phi^*\ran
=\lan \psi,\sigma\circ\phi^*\ran$  by \eqref{11} and thus
\begin{equation}\label{22}
\lan\psi\circ\phi,\sigma\ran
=\lan\phi,\psi^*\circ\sigma\ran
=\lan \psi,\sigma\circ\phi^*\ran.
\end{equation}

For a subset $K$ of $H(M_A,M_B)$ of Hermitian-preserving maps,
the corresponding \emph{dual cone} \markg{$K^\circ$ with respect to the bilinear pairing in \eqref{eq:maps_blinear_pair}} is the set defined as
$$
K^\circ=\{\psi\in H(M_A,M_B): \lan \phi,\psi\ran\ge 0\ {\text{\rm for every}}\ \phi\in K\}.
$$
We note that $K^{\circ\circ}$ is the smallest closed convex cone in $H(M_A,M_B)$ containing $K$.
In particular, we have $K^{\circ\circ}=K$ if and only if $K$  is a closed convex cone.
For closed convex cones $K_1$ and $K_2$, we denote by $K_1\join K_2$ and $K_1\meet K_2$
the convex hull and the intersection of $K_1$ and $K_2$, respectively. Then we have
$$
(K_1\join K_2)^\circ=K_1^\circ\meet K_2^\circ\qquad\text{and}\qquad
(K_1\meet K_2)^\circ=K_1^\circ\join K_2^\circ.
$$
(See \cite{{eom-kye},{kye_ritsu},{han_kye_par-sep}} for further properties of cones
in a more general setting). For a subset $K$ of $H(M_A,M_B)$, we also define the set $K^*$ as
\[
K^*=\{\phi^*\in H(M_B,M_A):\phi\in K\}.
\]
By the identity in \eqref{11}, we have
\[
K^{*\circ}=K^{\circ *}.
\]
For subsets $K_0\subset H(M_A,M_B)$, $K_1\subset H(M_B,M_C)$ and $K_2\subset H(M_A,M_C)$,
the identities in \eqref{22} yield the following equivalences:
\begin{multline}\label{lemma}
K_1\circ K_0\subset K_2^\circ
\ \Longleftrightarrow\
K_1^*\circ K_2\subset K_0^\circ
\ \Longleftrightarrow\
K_2\circ K_0^*\subset K_1^\circ\\
\ \Longleftrightarrow\
K_1\subset(K_2\circ K_0^*)^\circ
\ \Longleftrightarrow\
K_0\subset(K_1^*\circ K_2)^\circ
\ \Longleftrightarrow\
K_2\subset(K_1\circ K_0)^\circ,
\end{multline}
 where we define
\[
K_1\circ K_0=\{\phi_1\circ\phi_0\in H(M_A,M_C):\phi_1\in K_1, \phi_0\in K_0\}.
\]

We denote by $\cp_A$ the convex cone of all completely positive linear maps of $M_A$ into itself and note that
$\cp_A^\circ=\cp_A^*=\cp_A$. For a given closed convex cone $K$ in $H(M_A,M_B)$,
we are interested in conditions on $K$ that are equivalent to the conditions $K\circ \cp_A\subset K$
and $\cp_A\circ K\subset K$, respectively.
To do this, we may plug $K_0=\cp_A$, $K_1=K$ and $K_2=K^\circ$ into \eqref{lemma} to obtain the following equivalences:
\begin{multline}\label{tttttcp}
K\circ \cp_A\subset K
\ \Longleftrightarrow\
K^*\circ K^\circ\subset \cp_A
\ \Longleftrightarrow\
K^\circ\circ \cp_A\subset K^\circ\\
\ \Longleftrightarrow\
K\subset(K^\circ\circ \cp_A)^\circ
\ \Longleftrightarrow\
\cp_A\subset(K^*\circ K^\circ)^\circ
\ \Longleftrightarrow\
K^\circ\subset(K\circ \cp_A)^\circ.
\end{multline}
For a subset $K$ of $H(M_A,M_B)$, we define the sets $K^\varogreaterthan$ and $K^\varolessthan$ as
\[
K^\varogreaterthan =(K\circ\cp_A)^\circ\qquad\text{and}\qquad
K^\varolessthan     = (\cp_B\circ K)^\circ.
\]
For an arbitrary convex cone $K$ of $H(M_A,M_B)$ we have $K\subset K\circ\cp_A$ and
$K\subset \cp_B\circ K$, since the identity maps $1_A$ and $1_B$ are contained
in $\cp_A$ and $\cp_B$, respectively. This implies the inclusions
\begin{equation}\label{comparison}
K^\varogreaterthan\subset K^\circ\qquad\text{and}\qquad
K^\varolessthan\subset K^\circ.
\end{equation}
We also have $K^{\circ\varogreaterthan}\subset K^{\circ\circ}$ and $K^{\circ\varolessthan}
\subset K^{\circ\circ}$ by applying \eqref{comparison} to $K^\circ$.
The above analysis is summarized in \eqref{leftdual} of the following proposition.
The equivalences in \eqref{rightdual} can be shown analogously by choosing $K_0=K$, $K_1=\cp_B$ and $K_2=K^\circ$ in \eqref{lemma}.

\begin{proposition}\label{cccc}
For a closed convex cone $K$ in $H(M_A,M_B)$,
we have the following equivalent conditions:
\begin{equation}\label{leftdual}
\begin{aligned}
K^\circ=K^\varogreaterthan
\ &\Longleftrightarrow\
K^\circ\subset K^\varogreaterthan
\ \Longleftrightarrow\
K\circ \cp_A\subset K
\ \Longleftrightarrow\
K^\circ\circ\cp_A\subset K^\circ\\
\ &\Longleftrightarrow\
K^*\circ K^\circ\subset \cp_A
\ \Longleftrightarrow\
K\subset K^{\circ\varogreaterthan}
\ \Longleftrightarrow\
K= K^{\circ\varogreaterthan}.
\end{aligned}
\end{equation}
We also have the following:
\begin{equation}\label{rightdual}
\begin{aligned}
K^\circ=K^\varolessthan
\ &\Longleftrightarrow\
K^\circ\subset K^\varolessthan
\ \Longleftrightarrow\
\cp_B\circ K\subset K
\ \Longleftrightarrow\
\cp_B\circ K^\circ\subset K^\circ\\
\ &\Longleftrightarrow\
K^\circ \circ K^*\subset \cp_B
\ \Longleftrightarrow\
K\subset K^{\circ\varolessthan}
\ \Longleftrightarrow\
K= K^{\circ\varolessthan}.
\end{aligned}
\end{equation}
\end{proposition}

We call a closed convex cone $K$ of positive linear
maps in $H(M_A,M_B)$ a {\sl left-mapping cone} if $\cp_B\circ
K\subset K$ holds and a {\sl right-mapping cone} if $K\circ
\cp_A\subset K$. For a closed convex cone $K$ of positive maps in
$H(M_A,M_B)$, we denote by $\slmp_K$ (respectively $\srmp_K$) the
smallest left- (respectively right-) mapping cone that contains $K$.

\begin{proposition}\label{strjjjj11}
For a closed convex cone $K$ of positive maps, we have the following:
\begin{enumerate}
\item[(i)]
$\srmp_K=\mcone^{\varogreaterthan\circ}$ and this closed convex cone is generated by $K\circ\cp_A$.
\item[(ii)]
$\slmp_K=\mcone^{\varolessthan\circ}$ and this closed convex cone is generated by $\cp_B\circ K$.
\end{enumerate}
\end{proposition}

\begin{proof}
To prove (i), we first note that the  convex hull $(K\circ \cp_A)^{\circ\circ}$ of $K\circ\cp_A$ is a right-mapping cone.
For any other right-mapping cone $L$ satisfying $K\subset L$, we have
$K\circ \cp_A\subset L\circ \cp_A\subset L$ and thus $(K\circ\cp_A)^{\circ\circ} \subset L$. Hence
$\srmp_K=(K\circ\cp_A)^{\circ\circ}=\mcone^{\varogreaterthan\circ}$ and, moreover,
the cone $\mcone^{\varogreaterthan\circ}$ is generated by $K\circ\cp_A$. The proof of statement (ii) is analogous.
\end{proof}

The convex cones $K^\varogreaterthan$ and $K^\varolessthan$ can be described
in terms of compositions of maps, as is shown in the following proposition.

\begin{proposition}\label{zzzz}
For a subset $K$ of $H(M_A,M_B)$, we have the following identities:
$$
\begin{aligned}
K^\varogreaterthan
&=\{\phi\in H(M_A,M_B): \psi^*\circ \phi\in\cp_A\ {\text{\rm for every}}\ \psi\in K\}\\
K^\varolessthan
&=\{\phi\in H(M_A,M_B): \phi\circ\psi^*\in\cp_B\ {\text{\rm for every}}\ \psi\in K\}.
\end{aligned}
$$
\end{proposition}

\begin{proof}
To prove the first identity, note for every map $\phi\in H(M_A,M_B)$ that
\begin{align*}
\phi\in K^\varogreaterthan
\ &\Longleftrightarrow\
\lan \phi,\psi\circ\sigma\ran\ge 0\ {\text{\rm for every}}\ \psi\in K,\ \sigma\in\cp_A\\
\ &\Longleftrightarrow\
\lan \psi^*\circ\phi,\sigma\ran\ge 0\ {\text{\rm for every}}\ \psi\in K,\ \sigma\in\cp_A\\
\ &\Longleftrightarrow\
\psi^*\circ\phi\in\cp_A\ {\text{\rm for every}}\ \psi\in K
\end{align*}
by the identities in \eqref{22}. The proof of the second statement is analogous.
\end{proof}

\begin{theorem}\label{cor-mp-right}
For a closed convex cone $K$ of positive maps in $H(M_A,M_B)$, the following are equivalent:
\begin{enumerate}
\item[(i)]
$K$ is a right-mapping cone.
\item[(ii)]
$K^*$ is a left-mapping cone.
\item[(iii)] For all $\phi\in L(M_A,M_B)$,
$\phi\in K^\circ$ if and only if $\psi^*\circ\phi\in\cp_A$ for every $\psi\in K$.
\item[(iv)] For all $\phi\in L(M_A,M_B)$,
$\phi\in K$ if and only if $\psi^*\circ\phi\in\cp_A$ for every $\psi\in K^\circ$.
\end{enumerate}
We also have the following equivalent conditions:
\begin{enumerate}
\item[(v)]
$K$ is a left-mapping cone.
\item[(vi)]
$K^*$ is a right-mapping cone.
\item[(vii)]For all $\phi\in L(M_A,M_B)$,
$\phi\in K^\circ$ if and only if $\phi\circ\psi^*\in\cp_B$ for every $\psi\in K$.
\item[(viii)] For all $\phi\in L(M_A,M_B)$,
$\phi\in K$ if and only if $\phi\circ\psi^*\in\cp_B$ for every $\psi\in K^\circ$.
\end{enumerate}
\end{theorem}

\begin{proof}
By Proposition \ref{cccc}, we see that $K$ is a right-mapping cone if and only if
$K^\circ=K^\rdual$, which holds if and only if $K=K^{\circ\rdual}$.
We note that $K^\circ=K^\rdual$ is equivalent to (iii) and
$K=K^{\circ\rdual}$ is equivalent to (iv) by Proposition \ref{zzzz}.
To prove the equivalence (i)$\Longleftrightarrow$(ii), we first note that a map $\phi$ is positive if and only if $\phi^*$ is positive.
We also note that  $K^{\rdual *}=K^{*\ldual}$  for
every closed convex cone $K$ by Proposition \ref{zzzz}.
If $K$ is a right-mapping cone then $K^{*\circ}=K^{\circ *}=K^{\rdual *}=K^{*\ldual}$, and so we see that $K^*$ is a left-mapping cone
by Proposition \ref{cccc}. Similarly, if $K$ is a left-mapping cone then $K^*$ is a right-mapping cone.
The equivalences of (v), (vi), (vii) and (viii) are analogous.
\end{proof}

 The properties listed in (iii), (iv), (vii) and (viii) of Theorem \ref {cor-mp-right}
have been considered in
various contexts in the literature \cite{{stormer08-3},{sko-laa},{stormer_scand_2012},{stormer18}}.
We have shown here that such properties actually characterize
one-sided mapping cones. (See also Corollary \ref{cor-mp} in the next
section.) It should be noted that the dual cone $K^\circ$ of a left/right-mapping cone $K$
may not be a left/right-mapping cone, even though, by the equivalences in \eqref{tttttcp},
it holds that $K\circ \cp_A\subset K$  if and only if
$K^\circ\circ \cp_A\subset K^\circ$. The dual of a left/right-mapping cone may contain a non-positive map,
as we will see in Example \ref{one-sided-dual} of the next section.

We now investigate the Choi matrix $\choi_{\psi\circ\phi}$ of compositions for maps $\phi\in L(M_A,M_B)$ and $\psi\in L(M_B,M_C)$
in terms of the Choi matrices $\choi_\phi$ and $\choi_\psi$. This will be useful for the examples considered later in the paper.
For matrices $a\in M_A$ and $c\in M_C$, we use the identity in \eqref{dual-id} to see that
\begin{align*}
\lan a\ot c, \choi_{\psi\circ\phi}\ran_{AC}
&=\lan c,\psi(\phi(\markg{a})))\ran_{C}\\
&=\lan \psi^*(c),\phi(\markg{a})\ran_{B}\\
&=\sum_{k,\ell}\lan \psi^*(c),e_{k\ell}\ran_B\lan e_{k\ell},\phi(\markg{a})\ran_B\\
&=\sum_{k,\ell}\lan c,\psi(e_{k\ell})\ran_C\lan\phi^*(e_{k\ell}),\markg{a}\ran_A\\
&=\sum_{k,\ell}\lan a\ot c,\markg{\phi^*(e_{k\ell})}\ot\psi(e_{k\ell})\ran_{AC},
\end{align*}
where $\{e_{k\ell}\}$ denote the matrix units in $M_B$. We therefore have that
\begin{equation}\label{schur}
\choi_{\psi\circ\phi}
=\sum_{k\ell} \markg{\phi^*(e_{k\ell})}\ot \psi(e_{k\ell})\in M_A\ot M_C.
\end{equation}
This is the block-wise summation of the block Schur product \cite{{horn},{chri-schur}}
\begin{equation*}
\markg{\choi_{\phi^*}}\square \choi_\psi\in M_B\ot (M_A\ot M_C)
\end{equation*}
of
$\markg{\choi_{\phi^*}}\in M_B\ot M_A$ and $\choi_\psi\in M_B\ot M_C$.
We may also note that the Choi matrix of the composition may be given by
\begin{equation}\label{eq:choi_composition_simple}
 \choi_{\psi\circ\phi} = \sum_{k,\ell} e_{k\ell}\otimes\psi(\phi(e_{k\ell})) = (1_A\otimes\psi)(\choi_\phi).
\end{equation}


\begin{example}\label{example}
\em Consider the map $\sigma$ on the algebra of $2\times 2$
matrices defined by
\[
\sigma=1_{M_2}+\ad_{e_{21}},
\]
where $1_{M_2}$ is the identity map on $M_2$ and, for a fixed matrix $a$,
the map $\ad_a$ is the map defined by $\ad_a (x)=a^*xa$. The Choi matrices of $\sigma$ and
$\sigma^*$ are given by
$$
\choi_\sigma=\left(\begin{matrix}
1&\cdot&\cdot&1\\
\cdot&\cdot&\cdot&\cdot\\
\cdot&\cdot&1&\cdot\\
1&\cdot&\cdot&1
\end{matrix}\right)
\qquad{\text{\rm and}}\qquad
\choi_{\sigma^*}=\left(\begin{matrix}
1&\cdot&\cdot&1\\
\cdot&1&\cdot&\cdot\\
\cdot&\cdot&\cdot&\cdot\\
1&\cdot&\cdot&1
\end{matrix}\right),
$$
respectively, where $\cdot$ denotes a zero. We denote by $\mathbb P_1$
the convex cone consisting of all positive linear maps. We have the
relations $\mathbb P_1\circ\cp\subset\mathbb P_1$ and
$\cp\circ\mathbb P_1\subset \mathbb P_1$, and thus $\mathbb P_1$
satisfies all the conditions in Proposition \ref{cccc}. Define the cone $K$ as
\begin{equation}\label{exam}
K=(\mathbb P_1^\circ\join \{\sigma\})^{\circ\circ},
\end{equation}
which is the convex hull of $\mathbb P_1^\circ\join \{\sigma\}$.
Note that $K^\circ =\mathbb{P}_1\meet \{\sigma\}^\circ=\{\phi\in\mathbb{P}_1:\lan \phi,\sigma\ran\geq0\}$ and that
\begin{equation}\label{kkkkk}
\begin{aligned}
K^\varogreaterthan&=\mathbb P_1\meet \{\sigma\}^\varogreaterthan =\{\phi\in\mathbb P_1: \sigma^*\circ\phi\in\cp\}\\
K^\varolessthan&=\mathbb P_1\meet \{\sigma\}^\varolessthan =\{\phi\in\mathbb P_1: \phi\circ\sigma^*\in\cp\}\\
\end{aligned}
\end{equation}
by Proposition \ref{zzzz}. We will show that $K^\varolessthan\nsubseteq K^\varogreaterthan$ and $K^\varogreaterthan\nsubseteq K^\varolessthan$.

Toward this goal, consider maps $\phi_{[a,b,c,d]}:M_2\to M_2$ with Choi matrix having the form
\[
\choi_{\phi_{[a,b,c,d]}}
=\left(\begin{matrix}
a &\cdot &\cdot &-1 \\
\cdot &b&\cdot&\cdot \\
\cdot &\cdot&c&\cdot \\
-1 &\cdot &\cdot &d
\end{matrix}\right)
\]
for fixed non-negative numbers $a,b,c$ and $d$. (This is a $(2\times 2)$-variant of the generalized Choi map considered in~\cite{cho-kye-lee}.)
For every rank-one projection $|\xi\ran\lan\xi|$ having the form
$|\xi\ran=(x,y)^\ttt$ for fixed constants $x,y\in\mathbb{C}$, the map $\phi_{[a,b,c,d]}$ sends $\lvert \xi\ran\lan\xi\rvert$ to
\[
\phi_{[a,b,c,d]}(|\xi\ran\lan\xi|)=
\left(\begin{matrix}
a\lvert x\rvert^2+c\lvert y\rvert^2 &-x\overline{y}\\ -\overline{x}y & b\lvert x\rvert^2+d\lvert y\rvert^2
\end{matrix}\right).
\]
Therefore $\phi_{[a,b,c,d]}\in\mathbb P_1$ if and only if the above matrix is positive
for every $x,y\in\mathbb C$, which holds if and only if
$\sqrt{ad}+\sqrt{bc}\ge 1$. We also have that $\lan \phi_{[a,b,c,d]},\sigma\ran = a+c+d-2$ and that
$$
\choi_{\sigma^*\circ\phi}= \left(\begin{matrix}
a &\cdot &\cdot &-1 \\
\cdot &a+b&\cdot&\cdot \\
\cdot &\cdot&c&\cdot \\
-1 &\cdot &\cdot &c+d
\end{matrix}\right)\qquad\text{and}\qquad
\choi_{\phi\circ\sigma^*}= \left(\begin{matrix}
a+c &\cdot &\cdot &-1 \\
\cdot &b+d&\cdot&\cdot \\
\cdot &\cdot&c&\cdot \\
-1 &\cdot &\cdot &d
\end{matrix}\right).
$$
The identities in \eqref{kkkkk} therefore yield the following equivalences:
\begin{itemize}
\item
$\phi_{[a,b,c,d]}\in K^\circ$ if and only if $\sqrt{ad}+\sqrt{bc}\ge 1$ and $a+c+d\ge 2$.
\item
$\phi_{[a,b,c,d]}\in K^\rdual$ if and only if $\sqrt{ad}+\sqrt{bc}\ge 1$ and $a(c+d)\ge 1$.
\item
$\phi_{[a,b,c,d]}\in K^\ldual$ if and only if $\sqrt{ad}+\sqrt{bc}\ge 1$ and $(a+c)d\ge 1$.
\end{itemize}
Finally, we may take $[a,b,c,d]=[\frac{1}{3},1,1,1]$ to conclude that
$K^\ldual \nsubseteq K^\rdual$ and take $[a,b,c,d]=[1,1,1,\frac{1}{3}]$ to conclude that $K^\rdual \nsubseteq K^\ldual$.
\ $\square$
\end{example}

\section{Mapping cones revisited}\label{sec-mp}

Because the identity map is completely positive, we see that a closed convex cone~$K$ satisfies
$\cp_B\circ K\circ\cp_A\subset K$ if and only if it satisfies both $K\circ\cp_A\subset K$ and $\cp_B\circ K\subset K$.
Combining the equivalences in Proposition \ref{cccc} yields the following theorem.

\begin{theorem}\label{mp}
For a closed convex cone $K$ in $H(M_A,M_B)$, the following are equivalent:
\begin{enumerate}
\item[(i)]
$K^\circ = K^\varogreaterthan = K^\varolessthan$.
\item[(ii)]
$K = K^{\circ\varogreaterthan} = K^{\circ\varolessthan}$.
\item[(iii)]
$\cp_B\circ K\circ\cp_A\subset K$.
\item[(iv)]
$\cp_B\circ K^\circ \circ\cp_A\subset K^\circ$.
\item[(v)]
$K^*\circ K^\circ\subset \cp_A$ and $K^\circ\circ K^*\subset\cp_B$.
\end{enumerate}
\end{theorem}

Following \cite{stormer-dual}, we say that
a closed convex cone $\mcone\subset\mathbb P_1$ is a {\sl mapping cone} if
$\phi_1\circ\phi\circ\phi_2\in \mcone$ for all choices of maps $\phi\in \mcone$,
$\phi_1\in\cp_A$ and $\phi_2\in\cp_A$  (that is, if $K$ satisfies condition (iii) of Theorem \ref{mp}).
The following corollary now follows directly from Theorem \ref{cor-mp-right}.

\begin{corollary}\label{cor-mp}
For a closed convex cone $K\subset \mathbb P_1$, the following are equivalent:
\begin{enumerate}
\item[(i)]
$K$ is a mapping cone.
\item[(ii)]
$K$ is both a left- and right-mapping cone.
\item[(iii)] For all maps $\phi\in L(M_A,M_B)$,
$\phi\in K^\circ$ if and only if $\psi^*\circ\phi\in\cp_A$ for every $\psi\in K$
if and only if $\phi\circ\psi^*\in\cp_B$ for every $\psi\in K$.
\item[(iv)] For all maps $\phi\in L(M_A,M_B)$,
$\phi\in K$ if and only if $\psi^*\circ\phi\in\cp_A$ for every $\psi\in K^\circ$
if and only if $\phi\circ\psi^*\in\cp_B$ for every $\psi\in K^\circ$.
\end{enumerate}
\end{corollary}

The implications (i)$\Longrightarrow$(iii) and (i)$\Longrightarrow$(iv) of the above corollary are well known in various contexts
\cite{{stormer08-3},{sko-laa},{stormer_scand_2012},{stormer18}}.
Corollary \ref{cor-mp} tells us that the converses of these implications are also true.

The convex cone $\mathbb P_k$ consisting of all $k$-{\sl positive} linear maps is an example
of a mapping cone. Recall that a linear map $\phi\in L(M_A,M_B)$ is called $k$-{\sl positive}
if its ampliation $1_{M_k}\ot\phi$ with the $k\times k$ matrices
is positive.  The convex cone $\cp_{AB}$ of all completely positive maps in $L(M_A,M_B)$
coincides with $\mathbb P_{A\meet B}$, where $A\meet B$ denotes the minimum
of the dimensions of $\mathbb C^A$ and $\mathbb C^B$ (see Corollary 4.1.9 in \cite{stormer_book}).
The range of $\mathbb P_k$ under the Jamio\l kowski-Choi isomorphism is denoted by $\blockpos_k$.
Matrices in $\blockpos_k$ are called $k$-{\sl blockpositive}.
It was shown in \cite{eom-kye} that a map $\phi$ is $k$-positive if and only if
\[
\lan \choi_\phi, |\xi\ran\lan\xi|\ran_{AB}\ge 0\quad {\text{\rm for every }|\xi\ran
\in\mathbb{C}^A\otimes\mathbb{C}^B\text{ with Schmidt rank }\le k.}
\]
We denote by ${\mathcal S}_k$ the convex cone in $M_A\ot M_B$ that is
generated by all matrices of the form $|\xi\ran\lan\xi|$ for vectors $|\xi\ran\in\mathbb C^A\ot\mathbb C^B$
whose Schmidt rank is less than or equal
to $k$. For a matrix $X\in\pos_{AB}$, the smallest $k$ such that $X\in\mathcal{S}_k$
is called the \emph{Schmidt number} of $X$, where ${\mathcal P}_{AB}$ denotes
the set of of all positive matrices in $M_A\ot M_B$. The corresponding convex cone
in $L(M_A,M_B)$ of all maps whose Choi matrices have Schmidt number at most $k$
will be denoted by $\superpos_k$. Elements of $\superpos_k$ are said to be
$k$-{\sl superpositive} \cite{{ando-04},{ssz}}. Note that the convex cones $\mathbb P_k$ and $\superpos_k$ are
dual to each other and that ${\mathcal S}_{A\meet B}=\blockpos_{A\meet B}={\mathcal P}_{AB}$.
This duality is summarized by the following diagram, where JC denotes the Jamio\l{}kowsi-Choi isomorphism:
\begin{equation}\label{exam-conexx}
\begin{matrix}
L(M_A,M_B): & \superpos_1  &\subset &\superpos_k &\subset &\cp_{AB} &\subset &\mathbb P_k &\subset &\mathbb P_1\\
\\
\phantom{\jc}\downarrow\jc  &\downarrow  &&\downarrow  &&\downarrow  &&\downarrow &&\downarrow\\
\\
M_A\ot M_B: & {\mathcal S}_1  &\subset &{\mathcal S}_k &\subset &\pos_{AB} &\subset &\blockpos_k &\subset &\blockpos_1\\
\end{matrix}
\end{equation}

The convex cone ${\mathcal S}_1$ in the tensor product $M_A\ot M_B$ and its corresponding
convex cone $\superpos_1$ in $L(M_A,M_B)$ play crucial roles in quantum information theory.
Recall that a \emph{state} is a positive unital linear functional. Every state on the matrix algebra $M_A\ot M_B$
corresponds to a density matrix $\varrho\in M_A\ot M_B$ by the mapping $x\mapsto \lan x,\varrho\ran_{AB}$.
In this sense, we may identify a state with its corresponding density matrix. A density matrix in
$M_A\ot M_B$ is called {\sl separable} if it belongs to ${\mathcal S}_1$
and {\sl entangled} if it is not separable. By the duality between $\mathbb P_1$ and $\superpos_1$,
we see that a state $\varrho$ is entangled if and only if there exists $\phi\in\mathbb P_1$ such that $\lan\choi_\phi,\varrho\ran<0$
\cite{horo-1}.
Maps in $\superpos_1$ are also called {\sl entanglement-breaking} \cite{{shor},{hsrus},{ruskai_EB},{holevo_BE}}.
It is easy to see that the convex cone $\superpos_k$ is generated by maps of the form
$\ad_V$ for matrices $V$ with $\operatorname{rank}(V)\le k$.

The following results show that the condition $K\subset\mathbb P_1$
in the definition of a mapping cone may be weakened
significantly. Note that for any matrices $x\in M_A$ and $y\in M_B$, the mapping
$\psi\in L(M_A,M_B)$ defined as $\psi(z)=\lan  x,z\ran_A \, y$ for every $z\in M_A$ has Choi matrix equal to $\choi_\psi = x\otimes y$.

\begin{lemma}\label{MC----dual-lemma}
Let $\mcone$ be a nonzero mapping cone in $H(M_A,M_B)$ (i.e., a closed convex cone that satisfies
$\cp_B\circ \mcone\circ\cp_A\subset \mcone$) and suppose there exists a map $\phi\in\mcone$ such that $\tr(\choi_\phi)>0$.
Then we have
$\superpos_1\subset K$.
\end{lemma}

\begin{proof}
Let $x\in\pos_A$ and $y\in\pos_B$ and define maps $\psi\in L(M_A,M_A)$ and $\sigma\in L(M_B,M_B)$
whose Choi matrices are $\choi_\psi = x\otimes I_A$ and $\choi_\sigma = I_B\otimes y$ respectively.
These maps are completely positive and thus $\sigma\circ\phi\circ\psi\in\mcone$ by assumption.
By \eqref{eq:choi_composition_simple}, we have that
\[
 \choi_{\sigma\circ\phi\circ\psi} = \bigl(1_A\otimes (\sigma\circ\phi)\bigr)(\choi_{\psi})
 = x\otimes \sigma(\phi(I_A)) = \tr(\phi(I_A))\, x\otimes y = \tr(\choi_\phi)\, x\otimes y.
\]
The desired result now follows from the fact that $\pos_A\ot\pos_B$ generates the cone ${\mathcal S}_1$.
\end{proof}

\begin{proposition}\label{MC----dual}
Let $\mcone$ be a proper nonzero closed convex cone in  $H(M_A,M_B)$ that satisfies the condition
$\cp_B\circ \mcone\circ\cp_A\subset \mcone$. The following are equivalent:
\begin{enumerate}
\item[(i)]
$\mcone\subset\mathbb P_1$.
\item[(ii)]
there exists $\phi\in\mcone$ with $\tr(\choi_\phi)>0$.
\item[(iii)]
$\superpos_1\subset\mcone$.
\end{enumerate}
\end{proposition}

\begin{proof}
It is clear that (i) implies (ii), as every nonzero map $\phi\in\mathbb P_1$
satisfies $\tr(\choi_\phi)>0$.  The implication (ii)$\Rightarrow$(iii) follows from Lemma \ref{MC----dual-lemma}.
Suppose now that (iii) holds. Note that $K^\circ$ is nonzero by the assumption that $K$ is
proper and that (iii) is equivalent to $K^\circ\subset \mathbb P_1$, so there exists a map $\phi\in K^\circ$ satisfying $\tr(\choi_\phi)>0$.
By Theorem \ref{mp}, we may apply Lemma \ref{MC----dual-lemma} to the convex cone $K^\circ$
to see that $\superpos_1\subset K^\circ$, which implies that $K\subset \mathbb P_1$.
\end{proof}

It is known that every mapping cone $K$ satisfies $\superpos_1\subset K\subset \mathbb P_1$ (see Lemma 5.1.5 in \cite{stormer_book}).
From the equivalence (i)$\Longleftrightarrow$(iii) of Proposition \ref{MC----dual}
together with (iii)$\Longleftrightarrow$(iv) of Theorem \ref{mp}, we
recover the well-known fact \cite{{stormer_scand_2011},{sko-laa}} that
$K$ is a mapping cone if and only if $K^\circ$ is a mapping cone (see  Theorem 6.1.3 in \cite{stormer_book}).
In particular, the convex cone $\superpos_k$ is also a mapping cone.
The condition $\cp_B\circ \mcone\circ\cp_A\subset \mcone$ in Proposition \ref{MC----dual}
cannot be replaced by the weaker condition $\cp_B\circ K\subset K$, as is shown in the following example.

\begin{example}\label{one-sided-dual}
\em
For any fixed positive map $\sigma:M_A\to M_B$, the set $K=\cp_B\circ\{\sigma\}$ is a convex cone.
Moreover, it is clear that $K$ is a left-mapping cone and thus $\slmp_K=K$.
Consider now the map $\sigma$ whose Choi matrix is $e_{11}\ot e_{11}\in M_A\ot M_B$.
For an arbitrary $\phi\in L(M_A,M_B)$, we have
$$
\choi_{\phi\circ\sigma}=e_{11}\ot \phi(e_{11})\in M_A\ot M_B,
$$
by \eqref{eq:choi_composition_simple}. If the dimension of $\mathbb{C}^A$ is greater than 1,
it is clear that $\superpos_1\nsubset K$ and so we have that $K^\circ\nsubset \mathbb P_1$
even though $K\subset \mathbb P_1$.
\end{example}

If $\mcone$ is a mapping cone in $L(M_A,M_B)$ then it is easy to see that $\{\phi\circ\ttt:\phi\in \mcone\}$ is also a mapping cone,
where $\ttt$ denotes the transpose map.
We therefore obtain the following further examples of mapping cones defined by
\begin{equation}\label{co-exam}
\mathbb P^k:=\{\phi\circ\ttt: \phi\in\mathbb P_k\}\qquad\text{and}\qquad
\superpos^k:=\{\phi\circ\ttt: \phi\in\superpos_k\}.
\end{equation}
We have that $\mathbb P^{A\meet B}=\superpos^{A\meet B}$ and this cone will
be denoted by $\ccp$, whose elements are called {\sl completely copositive maps}.
It is also clear that the convex hull $\mcone_1\join\mcone_2$ and the intersection $\mcone_1\meet\mcone_2$
are mapping cones whenever $\mcone_1$ and $\mcone_2$ are mapping cones.
In particular, the mapping cones
$$
\dec:=\cp \join\ccp\qquad\text{and}\qquad \pptmap:=\cp\meet\ccp
$$
play important roles in the theory of quantum information. Elements of $\dec$ are called {\sl decomposable} positive maps
\cite{{stormer},{woronowicz},{stormer80}}.
It is known that $\mathbb P_1=\dec$ if and only if $(\dim \mathbb C^A,\dim \mathbb C^B)$ is $(2,2)$, $(2,3)$ or $(3,2)$
\cite{{stormer},{woronowicz},{choi75}}.
For a given $X\in M_A\ot M_B$, the matrix $X^\Gamma=(1_A\ot\ttt)(X)$ is called the {\sl partial transpose} of $X$.
For a map $\phi\in L(M_A,M_B)$, we see that
$$
\choi_{\phi\circ\ttt}=\sum_{ij}e^A_{ij}\ot\phi(e^A_{ji})=(\choi_\phi)^\Gamma.
$$
We therefore have that $\phi\in\pptmap$ if and only if both $\choi_\phi$ and $\choi_\phi^\Gamma$
are positive.
Such matrices are called \emph{positive partial transpose} (PPT).
Because ${\mathcal S}_1^\Gamma={\mathcal S}_1$, we have that $\varrho$ is PPT
for all separable states $\varrho\in\mathcal{S}_1$ \cite{{choi-ppt},{peres}}.
This is precisely the dual of the statement that $\dec\subset\mathbb P_1$.
A map $\phi\in L(M_A,M_B)$ is called a {\sl PPT} map if $\phi\in\pptmap$ (or equivalently $\choi_\phi$ is  PPT).
Therefore, a map is PPT if and only if
it is both completely positive and completely copositive.

Consider now the lattice generated by the mapping cones listed in \eqref{exam-conexx} and~\eqref{co-exam}
with respect to the following two operations: the convex hull $K_1\join K_2$ and the intersection $K_1\meet K_2$
of closed convex cones $K_1$ and $K_2$.
Mapping cones belonging to this lattice are said to be {\sl typical}~\cite{sko-laa}.
If $M_A$ is the set of $2\times 2$ matrices, this lattice of inclusions may be drawn
as one of the following diagrams  (depending on the dimension of $M_B$):
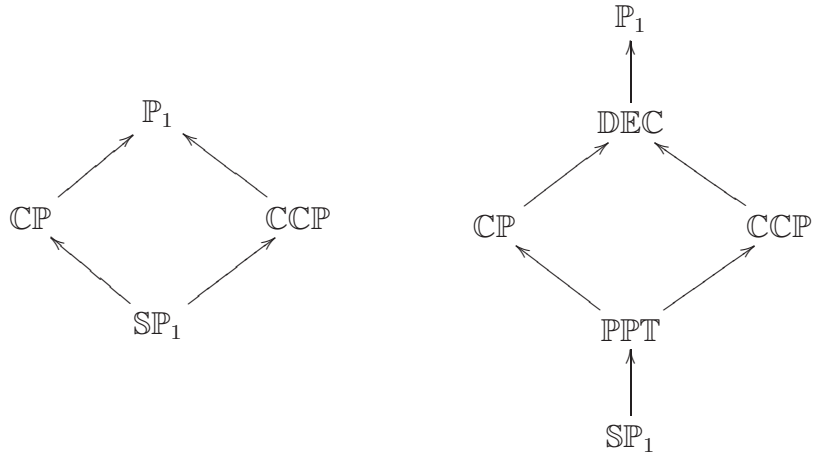
\[
\xymatrix{
\\
&\mathbb P_1& \\
\cp \ar[ur] && \ccp\ar[ul]\\
&\ar[ul]\superpos_1 \ar[ur]
}\qquad\qquad
\xymatrix{
&\mathbb P_1 & \\
&\dec \ar[u]& \\
\cp \ar[ur] && \ccp\ar[ul]\\
&\ar[ul]\pptmap \ar[ur]\\
&\superpos_1 \ar[u]&
}
\]
If the dimension of $\mathbb{C}^B$ is 2 or 3 then we have the diagram on the left.
If $\operatorname{dim}(\mathbb{C}^B)\geq4$ then we have the diagram on the right.
It is known that there exist mapping cones which are not typical \cite{johnston_untypical}.

\begin{example}\em
Consider again convex cone $K=(\mathbb{P}_1^\circ\meet \{\sigma\})^{\circ\circ}$
from Example \ref{example}. We will show that $K\subsetneqq K^{\ldual\ldual}$.
Toward this goal, consider the map $\tau:M_2\rightarrow M_2$ whose Choi matrix is given by
\[
\choi_\tau=
\left(\begin{matrix}
\cdot &\cdot&\cdot&\cdot\\
\cdot&1&1&\cdot\\
\cdot&1&1&\cdot\\
\cdot &\cdot&\cdot& 1
\end{matrix}\right).
\]
Note that we have $\tau^*=\tau$. For a map $\phi:M_2\rightarrow M_2$ with Choi matrix
given by $\choi_\phi=\left(\begin{matrix}x&y\\z&w\end{matrix}\right)$ for matrices $x,y,z,w\in M_2$,
the Choi matrices of the compositions $\phi\circ\tau^*$ and $\phi\circ\sigma^*$ are
$\choi_{\phi\circ\tau^*}=\left(\begin{matrix}w&z\\y&x+w\end{matrix}\right)$ and
$\choi_{\phi\circ\sigma^*}=\left(\begin{matrix} x+w &y\\ z& w\end{matrix}\right)$ respectively,
where $\sigma$ is the map from Example \ref{example}.
If it holds that $\phi\in K^\ldual$ then
$\choi_{\phi\circ\sigma^*}$ is positive and it is clear that
$\choi_{\phi\circ\tau^*}$ is also positive.
It follows that that $\choi_{\tau\circ\phi^*}$ is positive for every
$\phi\in K^\ldual$, from which we conclude that $\tau\in K^{\ldual\ldual}$.
Now, toward a contradiction, suppose that $\tau\in K$. There must exist a map $\phi\in\mathbb{SP}_1$
and a number $\lambda\geq0$ such that $\tau=\phi+\lambda\sigma$.
Comparing the Choi matrices $\choi_\tau$ and $\choi_\sigma$, we see that we must have $\lambda=0$
and thus $\tau\in\superpos_1$. However, it is clear $\choi_\tau$ is entangled as it is not PPT,
and thus $\tau\not\in\mathbb{SP}_1$. This is in contradiction to the assumption that $\tau\in K$.
We therefore conclude that $K\subsetneqq K^{\ldual\ldual}$.
\ $\square$
\end{example}

Let $K$ be an arbitrary closed convex cone of positive maps. Recall that the smallest left- (respectively right-) mapping cone
$\slmp_K$ (respectively $\srmp_K$) containing $K$
is given by $(\cp\circ K)^{\circ\circ}=K^{\ldual\circ}$ (respectively
$(K\circ\cp)^{\circ\circ}=K^{\rdual\circ}$).
By the same argument as in Proposition \ref{strjjjj11}, we see that the smallest mapping cone
$\smp_K$ containing $K$ is given by $\smp_K=(\cp\circ K\circ\cp)^{\circ\circ}$.
Moreover, it is clear that $(\cp\circ K)\cup (K\circ\cp)\subset \cp\circ K\circ\cp$.
We therefore obtain the following lattice of (not necessarily strict) inclusions:
\begin{equation}\label{examdiagramcc}
\xymatrix{
&  &\slmp_K \ar[rd]\\
\mcone \ar[r] &\slmp_K\meet \srmp_K \ar[ur]\ar[dr]
       &&\slmp_K\join \srmp_K \ar[r] &\smp_K\\
&&\srmp_K\ar[ru]
}
\end{equation}
Even though a closed convex cone $K$
is a mapping cone if and only if it is both a left- and right-mapping cone, we shall see
in the following example that $\slmp_K\join \srmp_K$ need not coincide with~$\smp_K$.
In particular, we will see that every inclusion in \eqref{examdiagramcc} is strict for the convex cone
$K$ from Example \ref{example}. For this cone, $\srmp_K$ is a right-mapping cone which is not a left-mapping cone,
as it holds that $\slmp_{\srmp_K}=\smp_K$ in general.

\begin{example}
\em
Consider again the closed convex cone $K=(\superpos_1\join\{\sigma\})^{\circ\circ}$,
where $\sigma$ is the map defined in Example \ref{example}. Taking the duals of
the convex cones in \eqref{examdiagramcc} yields the following chain of inclusions:
\begin{equation}\label{jkluy}
\xymatrix{
&  &K^\ldual \ar[rd]\\
(\cp\circ K\circ \cp)^{\circ} \ar[r] & K^\ldual\meet K^\rdual \ar[ur]\ar[dr]
       &&K^\ldual\join K^\rdual \ar[r] &K^\circ\\
&&K^\rdual\ar[ru]
}
\end{equation}
We will show that every inclusion in the above diagram is strict. It is clear
that every inclusion in the diamond part of the lattice is strict, since
$K^\ldual \nsubseteq K^\rdual$ and $K^\rdual \nsubseteq K^\ldual$ by Example~\ref{example}.
To see that the first inclusion is strict, consider the map $\phi_\alpha$ defined for a fixed positive number $\alpha>0$ by
\[
\phi_\alpha\left(\begin{pmatrix}
           x_{11} & x_{12}\\ x_{21} & x_{22}
          \end{pmatrix}
\right) = \left(\begin{matrix} \alpha x_{11}&x_{12}\\x_{21}&x_{22}/\alpha\end{matrix}\right).
\]
For $\alpha>0$ and $\beta>0$, we have that
\begin{equation}\label{kjhgfds}
\choi_{\phi_\beta\circ\sigma\circ\phi_\alpha}=
\left(\begin{matrix}
\alpha\beta &\cdot &\cdot&1\\
\cdot&\cdot&\cdot&\cdot\\
\cdot&\cdot&\beta/\alpha &\cdot\\
1&\cdot&\cdot&1/\alpha\beta
\end{matrix}\right)
\in \choi_{\cp\circ K\circ\cp}.
\end{equation}
Consider now the map having the form $\psi = \phi_{[a,b,c,d]}$ as defined in
Example \ref{example}, where we choose $a=b=c=d=\frac 1{\sqrt 2}$.
It is clear that $\psi\in K^\ldual\meet K^\rdual$. Toward a contradiction,
suppose that $\psi\in (\cp\circ K\circ \cp)^{\circ}$. It must be the case that
\[
\lan \psi,\phi_\beta\circ\sigma\circ\phi_\alpha\ran = \frac 1{\sqrt 2}\left(\alpha\beta+\frac\beta\alpha+\frac 1{\alpha\beta}\right)\ge 2
\]
for every $\alpha,\beta>0$ by \eqref{kjhgfds}. However, the above inequality fails to hold for $\alpha=2$ and $\beta=\frac 1{\sqrt 5}$.
It follows that $\psi\notin (\cp\circ K\circ \cp)^{\circ}$ and thus $(\cp\circ K\circ \cp)^{\circ} \subsetneqq K^\ldual\meet K^\rdual$.

Now consider maps of the form $\phi_{[a]}:=\phi_{[a,a,a,a]}$. We see that
$\phi_{[a]}\in\mathbb P_1$ if and only if $a\ge\frac 12$ and that $\phi_{[a]}\in K^\circ$ if and only if
$a\ge \frac 23$. On the other hand, we have $\phi_{[a]}\in K^\ldual$ if and only if
$\phi_{[a]}\in K^\rdual$ if and only if $a\ge \frac 1{\sqrt 2}$. It follows that the last inclusion
in \eqref{jkluy} is also strict.
Finally, we note that the mapping cone $\smp_K$ is not typical (see Theorem 18 of \cite{johnston_untypical}).
\ $\square$
\end{example}

\section{Duality through ampliation}\label{sec-amp}

In this section, we discuss relationships between one-sided mapping cones and ampliation maps $1_A\ot \phi$ and $\phi\ot 1_B$.
We will see that many dual objects---such as $K^\ldual$ and $K^\rdual$---can be described in terms of ampliation.
This allows us to recover many results in quantum information theory---such as
separability criteria through ampliation of positive maps \cite{horo-1} as well as characterizations of
entanglement-breaking maps \cite{hsrus} and Schmidt number \cite{terhal-sghmidt}---in a single framework.
We also recover some characterizations of decomposable maps due to the third author \cite{stormer82}
and properties of $k$-positive maps due to the second author \cite{eom-kye}. We stress that the above-mentioned characterizations hold
if and only if the involved convex cones are one-sided mapping cones.
For a convex cone $K$ in $H(M_A\ot M_B)$, we denote by $\choi_K$ the convex cone in the tensor product
$M_A\ot M_B$ defined as $\choi_K=\{\choi_\phi:\phi\in K\}$
.

\begin{proposition}\label{bbbbbbb}
For a closed convex cone $K$ in $H(M_A,M_B)$ and a map $\phi\in H(M_A,M_B)$, the following are equivalent:
\begin{enumerate}
\item[(i)]
$\phi\in K^\rdual$.
\item[(ii)]
$(1_A\ot\phi)(\choi_\psi)\in \choi_{K^\circ}$ for every $\psi\in\cp_{A}$.
\item[(iii)]
$(1_B\ot \psi)(\choi_{\phi^*})\in \choi_{K^{*\circ}}$ for every $\psi\in\cp_{A}$.
\item[(iv)]
$(1_A\ot\sigma^*)(\choi_\phi)\in \pos_{AA}$ for every $\sigma\in K$.
\item[(v)]
$(1_A\ot \phi^*)(\choi_\sigma)\in \pos_{AA}$ for every $\sigma \in K$.
\end{enumerate}
\end{proposition}

It is useful to remember the domains and the ranges of the ampliation maps in Proposition \ref{bbbbbbb}.
We have
$$
\begin{aligned}
1_A\ot\phi&:M_A\ot M_A\to M_A\ot M_B\\
1_B\ot \psi&: M_B\ot M_A\to M_B\ot M_A\\
1_A\ot\sigma^*&:M_A\ot M_B\to M_A\ot M_A\\
1_A\ot \phi^*&:M_A\ot M_B\to M_A\ot M_A
\end{aligned}
$$
for $\sigma\in K$ and $\psi\in\cp_A$.

\medskip
\begin{proof}
Note that statements (i) and (ii) are equivalent to the following two statements, respectively:
\begin{enumerate}
\item[(i${}^\prime$)]
$\lan \phi,\sigma\circ\psi^*\ran\ge 0$ for every $\psi\in\cp_{A}$ and $\sigma\in K$.
\item[(ii${}^\prime$)]
$\lan (1_A\ot\phi)(\choi_\psi),\choi_\sigma\ran_{AB}\ge 0$ for every $\psi\in\cp_{A}$ and $\sigma\in K$.
\end{enumerate}
From \eqref{eq:choi_composition_simple}, for every map $\psi:M_A\to M_A$ we have that
\begin{align*}
\lan \phi,\sigma\circ\psi^*\ran
&=\lan \phi\circ\psi,\sigma\ran
=\lan \choi_{\phi\circ\psi},\choi_\sigma\ran_{AB}
=\lan (1_A\ot\phi)(\choi_\psi),\choi_\sigma\ran_{AB},
\end{align*}
which proves the equivalence (i)$\Longleftrightarrow$(ii). Taking the flip of the above identity yields
$$
\lan \phi,\sigma\circ\psi^*\ran
=\lan \sigma^*,\psi^*\circ\phi^*\ran
=\lan\choi_{\sigma^*},(1_B\ot \psi^*)(\choi_{\phi^*})\ran_{BA},
$$
which proves the equivalence (i)$\Longleftrightarrow$(iii) since $\cp_A^*=\cp_A$.
On the other hand, we also have that
$$
\begin{aligned}
\lan \phi,\sigma\circ\psi^*\ran
&=\lan \sigma^*\circ\phi,\psi^*\ran
=\lan (1_A\ot\sigma^*)(\choi_\phi),\choi_{\psi^*}\ran_{AA}\\
&=\lan\psi,\phi^*\circ\sigma\ran
=\lan\choi_\psi,(1_A\ot\phi^*)(\choi_\sigma)\ran_{AA}.
\end{aligned}
$$
This completes the proof.
\end{proof}

Suppose that $M_A=M_B$.
For a fixed convex cone $K$ in $H(M_A,M_A)$, the third author \cite{stormer-dual} defined the set $P_K\subset M_A\otimes M_A$ as
$$
P_K=\{\varrho\in M_A\ot M_A: (1_A\ot\sigma)(\varrho)\ {\text {\rm is positive for every}}\ \sigma\in K\}.
$$
A map $\phi:M_A\to M_A$ is called $K$-{\sl positive} \cite{stormer-dual} if $\lan\varrho,\phi\ran\ge 0$ holds for every $\varrho\in K$.
By the equivalence of statements (i) and (iv) of Proposition \ref{bbbbbbb}, we see that $P_K=\choi_{K^{*\rdual}}$
and moreover that a map $\phi$ is $K$-positive if and only if $\phi\in K^{*\rdual\circ}$.
If it is the case that $K=K^*$, we see that a map $\phi:M_A\to M_A$ is $K$-positive if and only if $\phi\in K^{\rdual\circ}$,
which is equivalent to the condition that $\phi$ is the sum of maps of the form $\sigma\circ\psi$ for maps $\sigma\in K$ and $\psi\in\cp$
by statement (i) of Proposition \ref{strjjjj11}. This recovers a result from~\cite{stormer-dual},
where it was also shown that a map $\phi$ is $\cp$-positive if and only if $\phi\in\cp$.
This is a special case of the following characterization of right-mapping cones in terms of $K$-positivity,
which follows trivially form Proposition \ref{cccc}.

\begin{corollary}
Suppose that $K$ is a closed convex cone in $H(M_A,M_A)$ such that $K^*=K$. The following are equivalent:
\begin{enumerate}
\item[(i)]
For all maps $\phi:M_A\to M_A$, $\phi\in K$ if and only if $\phi$ is $K$-positive.
\item[(ii)]
$K\circ \cp_A\subset K$.
\end{enumerate}
\end{corollary}

We also have the following characterization of $K^\ldual$ by ampliation on the right.

\begin{proposition}\label{hhhh}
For a closed convex cone $K$ in $H(M_A,M_B)$ and a map $\phi\in H(M_A,M_B)$, the following are equivalent:
\begin{enumerate}
\item[(i)]
$\phi\in K^\ldual$.
\item[(ii)]
$(\phi^*\ot 1_B)(\choi_{\psi})\in \markg{\choi_{K^\circ}}$ for every $\psi\in\cp_B$.
\item[(iii)]
$(\psi\ot 1_A)(\markg{\choi_{\phi^*}})\in \markg{\choi_{K^{*\circ}}}$ for every $\psi\in\cp_B$.
\item[(iv)]
$(\sigma\ot 1_B)(\markg{\choi_\phi})\in \pos_{BB}$ for every $\sigma\in K$.
\item[(v)]
$(\phi\ot 1_B)(\markg{\choi_\sigma})\in \pos_{BB}$ for every $\sigma\in K$.
\end{enumerate}
\end{proposition}

\begin{proof}
For all maps $\sigma\in K$ and $\phi\in\cp_B$, we have the following maps:
\begin{align*}
\phi^*\ot 1_B&: M_B\ot M_B\to M_A\ot M_B \\
\psi\ot 1_A&:   M_B\ot M_A\to M_B\ot M_A \\
\sigma\ot 1_B&: M_A\ot M_B\to M_B\ot M_B \\
\phi\ot 1_B&:   M_A\ot M_B\to M_B\ot M_B.
\end{align*}
First recall from \eqref{eq:choi_composition_simple} that
\begin{equation}\label{eq:phisigmastar}
(1_B\ot\phi)(\choi_{\sigma^*})=
\choi_{\phi\circ\sigma^*}\in M_B\ot M_B.
\end{equation}
Recalling that $\choi_{\sigma\circ\phi^*}$ is the flip of $\markg{\choi_{\phi\circ\sigma^*}}$ from
\eqref{eq:ChoiPhiStarFlip}, taking the flip the identity in \eqref{eq:phisigmastar} yields
\[
(\phi\ot 1_B)(\markg{\choi_\sigma})
=\markg{\choi_{\sigma\circ\phi^*}}.
\]
We therefore obtain the following identities for all maps $\sigma\in K$ and $\phi\in\cp_B$:
\begin{align*}
\lan \phi,\psi\circ\sigma\ran
&=\lan\psi^*\circ\phi,\sigma\ran =\lan(\phi^*\ot 1_B)(\markg{\choi_{\psi^*}}),\markg{\choi_\sigma}\ran_{AB},\\
&=\lan\sigma^*,\phi^*\circ\psi\ran=\lan\markg{\choi_{\sigma^*}}(\psi^*\ot 1_A)(\markg{\choi_{\phi^*}})\ran_{BA}\\
&=\lan \phi\circ\sigma^*,\psi\ran
=\lan (\sigma\ot 1_B)(\markg{\choi_\phi}),\markg{\choi_\psi}\ran_{BB}\\
&=\lan\psi^*,\sigma\circ\phi^*\ran
=\lan \markg{\choi_{\psi^*}},(\phi\ot 1_B)(\markg{\choi_\sigma}) \ran_{BB}.
\end{align*}
Using the fact that $\cp_B=\cp_B^*$, applying an argument similar to the one in the proof of Proposition \ref{bbbbbbb} yields the desired conclusion.
\end{proof}

We now apply Proposition \ref{cccc} to Proposition \ref{bbbbbbb} and Proposition \ref{hhhh} to obtain the following
characterizations of one-sided mapping cones in terms of ampliation maps.

\begin{theorem}\label{mainth}
For a closed convex cone $K$ in $H(M_A,M_B)$, the following are equivalent:
\begin{enumerate}
\item[(i)]
$K\circ \cp_{AA}\subset K$.
\item[(ii)] For all maps $\phi$,
$\phi\in K^\circ$ if and only if
$(1_A\ot\phi)(\choi_\psi)\in \choi_{K^\circ}$ for every $\psi\in\cp_{A}$.
\item[(iii)]For all maps $\phi$,
$\phi\in K^\circ$ if and only if
$(1_B\ot \psi)(\choi_{\phi^*})\in \choi_{K^{*\circ}}$ for every $\psi\in\cp_{A}$.
\item[(iv)]For all maps $\phi$,
$\phi\in K^\circ$ if and only if
$(1_A\ot\sigma^*)(\choi_\phi)\in \pos_{AA}$ for every $\sigma\in K$.
\item[(v)]For all maps $\phi$,
$\phi\in K^\circ$ if and only if
$(1_A\ot \phi^*)(\choi_\sigma)\in \pos_{AA}$ for every $\sigma \in K$.
\end{enumerate}
We also have the following equivalence statements:
\begin{enumerate}
\item[(vi)]
$\cp_{BB}\circ K\subset K$.
\item[(vii)] For all maps $\phi$,
$\phi\in K^\circ$ if and only if
$(\phi^*\ot 1_B)(\markg{\choi_{\psi}})\in \markg{\choi_{K^\circ}}$ for every $\psi\in\cp_B$.
\item[(viii)] For all maps $\phi$,
$\phi\in K^\circ$ if and only if
$(\psi\ot 1_A)(\markg{\choi_{\phi^*}})\in \markg{\choi_{K^{*\circ}}}$ for every $\psi\in\cp_B$.
\item[(ix)]For all maps $\phi$,
$\phi\in K^\circ$ if and only if
$(\sigma\ot 1_B)(\markg{\choi_\phi})\in \pos_{BB}$ for every $\sigma\in K$.
\item[(x)] For all maps $\phi$,
$\phi\in K^\circ$ if and only if
$(\phi\ot 1_B)(\markg{\choi_\sigma})\in \pos_{BB}$ for every $\sigma\in K$.
\end{enumerate}
\end{theorem}

If $K$ is a mapping cone then statements (i) and (vi) in Theorem \ref{mainth} are trivially true.
We therefore have the following characterization of mapping cones.

\begin{corollary}\label{maincor}
Let $K$ be a mapping cone in $H(M_A,M_B)$. For a map $\phi\in L(M_A,M_B)$,
the following statements are equivalent to the statement that $\phi\in K$:
\begin{enumerate}
\item[(i)]
$1_A\ot\phi$ sends $\pos_{AA}$ into $\choi_K$,
\item[(ii)]
$1_B\ot \psi$ sends $\choi_{\phi^*}$ into $\choi_{K^*}$ for every $\psi\in\cp_{A}$,
\item[(iii)]
$1_A\ot\sigma^*$ sends $\choi_\phi$ into $\pos_{AA}$ for every $\sigma\in K^\circ$,
\item[(iv)]
$1_A\ot \phi^*$ sends $\choi_{K^\circ}$ into $\pos_{AA}$,
\item[(v)]
$\phi^*\ot 1_B$ sends $\pos_{BB}$ into $\markg{\choi_K}$,
\item[(vi)]
$\psi\ot 1_A$ sends $\markg{\choi_{\phi^*}}$ into $\markg{\choi_{K^*}}$ for every $\psi\in\cp_B$,
\item[(vii)]
$\sigma\ot 1_B$ sends $\markg{\choi_\phi}$ into $\pos_{BB}$ for every $\sigma\in K^\circ$,
\item[(viii)]
$\phi\ot 1_B$ sends $\markg{\choi_{K^\circ}}$ into $\pos_{BB}$.
\end{enumerate}
\end{corollary}

Applying statement (viii) to the cone $K=\dec$ allows us to recover the result in~\cite{stormer82},
which states that a map $\phi\in H(M_A,M_B)$ is decomposable if and only if
$\phi\ot 1_B$ sends PPT matrices in $M_A\ot M_B$ into the cone of positive matrices $\mathcal{P}_{BB}$.
Applying statement (viii) to the cone $K=\mathbb P_k$ shows that a map $\phi$ is $k$-positive if and only if
$\phi\ot 1_B$ sends every matrix with Schmidt number at most $k$ to a positive matrix, which recovers
a result in \cite{eom-kye}. The notion of $k$-positivity can be also characterized in terms
of the left-side ampliation as shown by the following corollary.

\begin{corollary}
Let $\phi\in L(M_A,M_B)$ be a map. The condition that $\phi$ is $k$-positive is equivalent to each of the following statements:
\begin{enumerate}
\item[(i)]
$1_A\ot\phi:M_A\ot M_A\to M_A\ot M_B$ sends ${\mathcal P}_{AA}$ into $\blockpos_k$,
\item[(ii)]
$1_B\ot\psi:M_B\ot M_A\to M_B\ot M_A$ sends $\choi_{\phi^*}$ into $\blockpos_k$ for every $\psi\in\cp_A$,
\item[(iii)]
$1_A\ot\sigma^*:M_A\ot M_B\to M_A\ot M_A$ sends $\choi_\phi$ into ${\mathcal P}_{AA}$ for every $\sigma\in\superpos_k$,
\item[(iv)]
$1_A\ot \phi^*:M_A\ot M_B\to M_A\ot M_A$ sends ${\mathcal S}_k$ into ${\mathcal P}_{AA}$.
\end{enumerate}
\end{corollary}

Taking the cone $K=\superpos_1$ and applying Corollary \ref{maincor}, we see that a map $\phi$
is entanglement breaking (i.e., $\choi_\phi$ is separable)  if and only if
$1_A\ot\phi$ sends every state to a separable state. This recovers a result from \cite{hsrus}.
A similar result for $k$-superpositive maps, which can be found in \cite{cw-EB}, can be stated as follows.
A map $\phi$ is $k$-superpositive if and only if $\phi\ot 1_B$ sends every $k$-blockpositive matrix to
a positive matrix.

We may interpret Corollary \ref{maincor} in terms of Choi matrices $\choi_\phi$ instead of the map $\phi$ itself.
For example, for the convex cone $K=\superpos_1$, from statement (vii) we have that
$\varrho\in \choi_{\superpos_1}={\mathcal S_1}$ if and only if
$\sigma\ot 1_B$ sends $\varrho$ to a positive matrix for every $\sigma\in \mathbb P_1=\superpos_1^\circ$.
That is, we see that a state $\varrho$ is separable if and only if $(\sigma\ot 1_B)(\varrho)$ is positive
for every positive map $\sigma$ \cite{horo-1}. Similarly, we also have that $\varrho$ has Schmidt number at most $k$
if and only if $(\phi\ot 1_B)(\varrho)$ is positive for every $k$-positive map $\sigma$ \cite{terhal-sghmidt}.
We conclude this section by presenting the following further characterizations of separability
that are found by applying statements (ii), (iii), (vi) and (vii) of Corollary \ref{maincor}. (See, e.g., \cite{horo-1}.)

\begin{corollary}
Let $\varrho\in M_A\ot M_A$ be a state. The condition that is $\varrho$ separable is equivalent to each of the following statements:
\begin{enumerate}
\item[(i)]
$\psi\ot 1_B$ sends $\varrho$ into ${\mathcal S}_1$ for every $\psi\in\cp_{A}$.
\item[(ii)]
$1_A\ot\sigma^*$ sends $\varrho$ into $\pos_{AA}$ for every $\sigma\in \mathbb P_1$.
\item[(iii)]
$1_A\ot \psi$ sends $\varrho$ into ${\mathcal S}_1$ for every $\psi\in\cp_B$.
\item[(iv)]
$\sigma\ot 1_B$ sends $\varrho$ into $\pos_{BB}$ for every $\sigma\in \mathbb P_1$.
\end{enumerate}
\end{corollary}

\section{PPT-square conjecture}\label{sec-PPT}

The notion of positive partial transpose plays an important role
in quantum information theory, as evidenced by the PPT criterion
for separability (${\mathcal S}_1\subset\ppt$).
The following conjecture was recently proposed by Christandl in \cite{ppt}.

\begin{conjecture}\label{ppt}
If $\phi$ and $\psi$ are PPT maps in $M_A$ then $\psi\circ\phi$ is entanglement breaking.
\end{conjecture}

This conjecture is called the PPT-square conjecture.
In our notation, the conjecture can be stated as the following inclusion:
$$
\pptmap\circ\pptmap\subset\superpos_1.
$$
The conjecture has been supported by the following results.
If $\phi$ is a unital or trace preserving PPT map then $\displaystyle{\lim_{k\to\infty}}d(\phi^k, \superpos_1)\to 0$ \cite{kennedy17}.
If $\phi$ is a unital PPT map then there is a positive integer $n$ such that $\phi^n\in \superpos_1$ \cite{rahaman18}.
Moreover, the conjecture has been shown recently to be true in the case when $M_A$ is the set of $3\times 3$ matrices
\cite{{chen_yany_tang},{Christandl19}}.
See also \cite{{collins_PPT},{hanson_PPT},{Muller-Hermes18}} for related results.

Choosing the cones $K_0=K_1=\pptmap$ and $K_2=\mathbb P_1$, applying the equivalences in \eqref{lemma} yields the following equivalences:
$$
\pptmap\circ\pptmap\subset\superpos_1\
\Longleftrightarrow\
\pptmap \circ \mathbb P_1 \subset \dec\
\Longleftrightarrow\
\mathbb P_1 \circ \pptmap \subset \dec.
$$
The equivalence of the first and the third of the above statements was shown in~\cite{Christandl19}.
From the identity $\pptmap=\pptmap\circ\cp=\pptmap\circ\dec$, we also have the following equivalences:
\[
\pptmap\circ\pptmap\subset\superpos_1\
\Longleftrightarrow\
\pptmap\circ\cp\circ\pptmap\subset\superpos_1\
\Longleftrightarrow\
\pptmap\circ\dec\circ\pptmap\subset\superpos_1.
\]
By the identity in \eqref{22}, we have that
\[
\lan \phi_0\circ\phi_1\circ\phi_2,\phi_3\ran =\lan \phi_1, \phi_0^*\circ\phi_3\circ\phi_2^*\ran,
\]
for all maps $\phi_0\in K_0$, $\phi_1\in K_1$, $\phi_2\in K_2$, and $\phi_3\in K_3$, which yields the equivalence
$$
K_0\circ K_1\circ K_2\subset K_3^\circ\
\Longleftrightarrow\
K_0^*\circ K_3\circ K_2^*\subset K_1^\circ
$$
for arbitrary closed convex cones $K_0$, $K_1$ $K_2$ and $K_3$.
This observation yields the equivalences
\begin{align*}
\pptmap \circ \cp \circ \pptmap \subset \superpos_1\
&\Longleftrightarrow\
\pptmap \circ \mathbb P_1 \circ \pptmap \subset \cp,\\
\pptmap \circ \dec \circ \pptmap \subset \superpos_1\
&\Longleftrightarrow\
\pptmap \circ \mathbb P_1 \circ \pptmap \subset \pptmap.
\end{align*}
The above equivalences may be summarized by the following theorem.
\begin{theorem}
The following statements are equivalent:
\begin{enumerate}
\item[(i)]
$\pptmap \circ \pptmap \subset \superpos_1$.
\item[(ii)]
$\pptmap \circ \mathbb P_1 \subset \dec$.
\item[(iii)]
$\mathbb P_1 \circ \pptmap \subset \dec$.
\item[(iv)]
$\pptmap \circ \cp \circ \pptmap \subset \superpos_1$.
\item[(v)]
$\pptmap \circ \dec \circ \pptmap \subset \superpos_1$.
\item[(vi)]
$\pptmap \circ \mathbb P_1 \circ \pptmap \subset \cp$.
\item[(vii)]
$\pptmap \circ \mathbb P_1 \circ \pptmap \subset \pptmap$.
\end{enumerate}
\end{theorem}

Using the identity in \eqref{schur}, we may also formulate the PPT-square conjecture in terms of block matrices as follows.

\begin{conjecture}
For all states $\varrho_1, \varrho_2\in M_A\ot M_A$, if both $\varrho_1$ and $\varrho_2$ are PPT then
the block-wise summation of the block Schur product $\varrho_1\square \varrho_2\in M_A\ot (M_A\ot M_A)$ is separable.
\end{conjecture}

\end{document}